\documentclass{article}
\usepackage{aism-2e}
\usepackage{amsmath,amssymb}
 
\newtheorem{theorem}{\sc Theorem}
\newtheorem{lemma}{\sc Lemma}

\newtheoremrm{definition}{\sc Definition}[thrmctr]
\newtheoremrm{assumption}{\sc Assumption}[thrmctr]
\newtheoremrm{condition}{\sc Condition}[thrmctr]
\newtheoremrm{remark}{\it Remark\/}                  
\newtheoremrm{example}{\it Example\/}
\newtheoremrm{case}{\it Case\/}
\newtheoremrmno{proof}{\sc Proof}
\newtheoremrmno{procedure}{\sc Procedure}
\newtheoremrmno{note}{\it Note\/}
\newtheoremrmno{fact}{\it Fact\/}
\makeatletter
\@addtoreset{corollary}{section}
\@addtoreset{proposition}{section}
\@addtoreset{assumption}{section}
\@addtoreset{condition}{section}
\frompage{1}\topage{12}

 \newcommand{\cC}{\mathcal{C}} 
  
 \newcommand{\cF}{\mathcal{F}}  
 \newcommand{\cG}{\mathcal{G}}

 \newcommand{\cN}{\mathcal{N}}
 \newcommand{\cR}{\mathcal{R}}
 \newcommand{\bC}{\mathbb{C}}
 \newcommand{\bR}{\mathbb{R}}
\newcommand{\bN}{\mathbb{N}}

\begin{document}
\sloppy

\title{
Integral Representations for Convolutions of \\ 
Non--Central Multivariate Gamma Distributions}

\author{
T. Royen
}
\affiliation{
Fachhochschule Bingen, University of applied sciences,\\
 Berlinstrasse 109, D--55411 Bingen, Germany\\
 E--mail: royen@fh-bingen.de
}


\abstract{
Three types of integral representations for the cumulative distribution functions of convolutions of $\Gamma_p(\alpha_k, \Sigma_k, \Delta_k)$--distributions with non--centrality matrices $\Delta_k$ are given by integration of products of simple complex functions over the $p$--cube $(-\pi,\pi]^p$. In particular, the joint distribution of the diagonal elements of a generalized quadratic form $XAX'$ with $n$ independent $\cN_p(\mu_k,\Sigma)$--distributed columns in $X_{p\times n}$ and a fixed $A \geq 0$ is obtained. For a single $\Gamma_p(\alpha, \Sigma,\Delta)$--cdf $(p-1)$--variate integrals over $(-\pi,\pi]^{p-1}$ are derived. The integrals are numerically more favourable than integrals obtained from the Fourier-- or Laplace inversion formula.} 

\date{}
\maketitle

\section{\large Introduction }
\label{section:1}
The following notations are used: $\sum_{(n)}$ stands for $\sum_{n_1+\ldots+n_p=n}$ with \mbox{$n_1,\ldots,n_p \in \bN_0$} and $\sum$ without any indices means $\sum ^{\infty}_{n=0}\sum_{(n)}$. The notation $D \geq 0$ is also used for non--symmetrical matrices $D_{p \times p}$ with only non--negative eigenvalues. The spectral norm of a $p \times p$--matrix $B$ is denoted by $\|B\|$, $I$ or $I_p$ is always an identity matrix and $\cC_p$ is the $p$--cube  $(-\pi,\pi]^p$.

The Laplace transform (L.t.) of a $p$--variate non--central $\Gamma_p(\alpha, \Sigma,\Delta)$--density with $\alpha >0$, $\Sigma >0$ and a non--centrality matrix $\Delta \geq 0$ was originally obtained from the L.t. of a non-central $W_p(2\alpha, \Sigma,\Delta)$--Wishart distribution (with an additional scale factor 2) and is given by 
\begin{eqnarray}
\label{eq:1}	
 \widehat{f}(t_1,\ldots,t_p;\alpha, \Sigma, \Delta)&= & |I_p + \Sigma T|^{-\alpha} \mathrm{etr} (-\Sigma T(I+\Sigma T)^{-1} \Delta),\\
T &=& \diag (t_1,\ldots,t_p) , t_1,\ldots,t_p \geq 0. \nonumber
\end{eqnarray}This function $\widehat{f}$ is generally the L.t. of the density of a real measure on $(0, \infty)^p$ which is not always a probability measure. The term "$\Gamma_p(\alpha, \Sigma,\Delta)$--distribution" is used here in this general sense. The exact set of values $\alpha$, leading to a probability density (pdf) $f(x_1,\ldots,x_p;\alpha, \Sigma, \Delta)$, depends on $\Sigma $ and presumedly on $\Delta$. To obtain a pdf, all positive integers $2\alpha$ (degrees of freedom) are admissible and all $2\alpha > p-1$. Moreover, in the central case all non--integer values $2\alpha > p-2\geq 0$ are allowed. For $p-2 < 2\alpha < p-1$ see Royen (1997). Furthermore all $\alpha > 0$ are admissible if $|I + \Sigma T|^{-1}$ is infinitely  divisible. Two characterizations of infinite divisibility of a $\Gamma_p(\alpha, \Sigma)$--distribution are found in Griffiths (1984) and Bapat (1989). Further conditions for admissible non--integer $2\alpha < p-2$ are given in Royen (1997), (2006).

Three integral representations by integration over $\cC_p$ are provided by theorem~2     
in section~4 for the functions 
\begin{eqnarray}
\label{eq:2}	
\lefteqn{F(x_1,\ldots,x_p; \alpha_1,\ldots,\alpha_n, \Sigma_1, \ldots,\Sigma_n, \Delta_1,\ldots,\Delta_n)}\\
&=& \int^{x_1}_{0}\ldots \int^{x_p}_{0} f(\xi_1,\ldots,\xi_p;\alpha_1,\ldots,\alpha_n, \Sigma_1,\ldots,\Sigma_n, \Delta_1,\ldots,\Delta_n)d\xi_1\ldots d\xi_p, \nonumber
\end{eqnarray}
where $f$ has the L.t.
\begin{eqnarray}
\label{eq:3}
\lefteqn{\prod^{n}_{k=1}	|I_p + \Sigma_kT|^{-\alpha_k} \mathrm{etr} (-\Sigma_k T(I + \Sigma_kT)^{-1} \Delta _k),}\\
&& \alpha_1,\ldots,\alpha_n > 0, \Sigma_1,\ldots,\Sigma_n > 0, \Delta_1,\ldots,\Delta_n \geq 0 . \nonumber
\end{eqnarray}
Thus, $F$ is not always the cumulative distribution function (cdf) of a probability measure. 

In particular let $X_{p\times n}$ be a $\cN_{p \times n}(M_{p \times n }, \Sigma_{p\times p} \otimes I_n)$--random matrix and $A_{n \times n} \geq 0$ of rank $q$ with $T'AT = \Lambda = \mathrm{diag} (\lambda_1,\ldots,\lambda_n)$, $\lambda_1\geq \ldots\geq \lambda_n$. Then the joint distribution of the diagonal elements of the generalized quadratic form $\frac 12 XAX'$ equals the distribution of the diagonal of $\frac12 Y\Lambda Y' $ with a $\cN_{p\times n} (MT, \Sigma_{p\times p} \otimes I_n)$--distributed  $Y = XT$. This is the distribution of a sum of $q$ independent $\Gamma_p(\frac12, \lambda_k\Sigma, \Delta_k = \frac12 \mu_k^* \mu_k^{*'} \Sigma^{-1})$--random vectors, where $\mu_k^* $ is the $k$--th column of $M^* = MT$. This joint distribution of $p$ quadratic forms of normal random vectors is comprised within theorem~2 as a special case with $\alpha_k = \frac12$, $\Sigma_k = \lambda_k \Sigma$, $k = 1,\ldots,q$. For methods under more general assumptions see also Blacher (2003). For a survey of univariate quadratic forms of normal random variables see chapter~4 in Mathai and Provost (1992). For several quadratic forms of skew elliptical distributions see B.Q. Fang (2005).

In Royen (1991), (1992) three different types of series expansions for the $\chi^2_p(2\alpha, \Sigma)$--cdf were derived from three different representations of the $\chi^2_p(2\alpha, \Sigma)$--L.t. which are extended to the general $\Gamma_p(\alpha,\Sigma,\Delta)$--L.t. in section~3 in a similar way as in Royen (1995). 

Some series expansions, closely related to the first two types, are already found in Khatri, Krishnaiah and Sen (1977). The third type was introduced because of its superior convergence properties.  The simple method to transform many series expansions into integrals over $\cC_p$ is explained in more detail in section~2 and summarized in theorem~1. The idea is as follows:

If $A(z_1,\ldots,z_p)$ and $B(z_1,\ldots,z_p)$ are analytical functions whose power series have the coefficients $a(m_1,\ldots,m_p)$ and $b(n_1,\ldots,n_p)$ and which are absolutely convergent for $\max |z_j| < r_A$ and $\max |z_j| < r_B$ respectively, where $r^{-1}_B < r_A$, then
\begin{eqnarray}
\label{eq:4}
\lefteqn{(2\pi)^{-p} \int_{\cC_p} A(y_1,\ldots,y_p) B(y^{-1}_{1},\ldots,y^{-1}_{p})d\varphi_1 \ldots d\varphi_p}\\
&& = \sum a(n_1,\ldots,n_p) b(n_1,\ldots,n_p)\nonumber
\end{eqnarray}
holds with $y_j = re^{i\varphi_j}$, $-\pi < \varphi_j \leq \pi$, $j = 1,\ldots,p$ and $r^{-1}_B < r < r_A$.

The integrals in (\ref{eq:4}) might be more economical than the series if the generating functions $A$ and $B$ are simple available functions and if the series are slowly convergent with very intricate coefficients. For non--central multivariate gamma distributions series expansions are practically not feasible.

The integral representations in theorem~2 of section~4 are of the type in (\ref{eq:4}). As long as no elementary density formulas are availale it should be a reasonable way to obtain the joint cdf by integration of elementary terms only over $\cC_p$ and not over $\bR^p$ as by the Fourier or Laplace inversion formula.  A single $\Gamma_p(\alpha,\Sigma,\Delta)$--cdf is represented by a $(p-1)$--variate integral over $\cC_{p-1}$ in section~5.

A totally different ${m+1\choose 2}$--variate integral representation of the $\Gamma_p(\alpha, \Sigma)$--cdf has been given recently by Royen (2006), which is based on $m$--factorial decompositions \mbox{$\sum^{-1}_{p\times p}  = D - BB'$,} where $D$ is a real or complex diagonal matrix minimizing the rank $m$ of $\Sigma^{-1} - D$. Approximations to a $\Gamma_p(\alpha, \Sigma)$--cdf are obtained by $m$--factorial approximations to $\Sigma$ with a low value of $m$. These approximations are improved further by successive correction terms.

\section{\large The method }
\label{section:2}
Theorem~1 in this section can be generalized in many ways, e.g. for Fourier transforms, but the version below is sufficient for the purpose of the underlying paper.

Let $\widehat{f}(t_1,\ldots,t_p)$, $t_1,\ldots,t_p \geq 0$, be a given L.t. of an unknown function $f(x_1, \ldots,x_p)$ with $f=0$ for $\min x_j < 0$. It is assumed that there are univariate L.t. $\widehat{g}_{j_0} (t)$ of some probability densities $g_{j_0}(x)$ on $(0,\infty)$ and further functions $h_j(t)$ with $|h_j(t)| \leq 1$, uniformly for $t \geq 0$, which enable a representation
\begin{eqnarray}
\label{eq:5}
\widehat{f}(t_1,\ldots,t_p) &=& \left( \prod^{p}_{j=1} \widehat{g}_{j_0}(t_j)\right) B\left(h_1(t_1),\ldots,h_p(t_p)\right)
\end{eqnarray}
with an analytical function $B(z_1,\ldots,z_p)$ whose power series expansion
\begin{eqnarray}
\label{eq:6}
\sum b(n_1,\ldots,n_p) \prod^{p}_{j=1} z^{n_j}_{j}
\end{eqnarray}
is absolutely convergent for $|z_1|,\ldots,|z_p| < r_B$ with a certain value $r_B > 1$.

Furthermore, the products $\widehat{g}_{j_0}(t) (h_j(t))^n$ are supposed to be the L.t. of continuous functions $g_{jn}(x)$, $x >0$, which satisfy the conditions 
\begin{eqnarray}
\label{eq:7}
&&|g_{jn}(x)|\leq  n^c k(x) \mbox{ with a constant $c$ and}\nonumber\\[-2ex]
&&\\[-2ex]
&&\int^{\infty}_{0}k(x)e^{-tx}dx < \infty \mbox{ for  all $ t > 0$ .}\nonumber
\end{eqnarray}
Hence, the generating functions (generators)
\begin{eqnarray}
\label{eq:8}
g_j(x,y) = \sum^{\infty}_{n=0} g_{jn}(x) y^n, \quad j=1,\ldots,p,
\end{eqnarray}
are defined for all $x > 0$ and $|y|<1$, and they have the L.t.
\begin{eqnarray}
\label{eq:9}
\widehat{g}_j(t,y) = \frac{\widehat{g}_{j_0}(t)}{1-yh_j(t)}, \quad t \geq 0 .
\end{eqnarray}

\begin{theorem}
\label{theorem:1}
Under the assumptions from (\ref{eq:6}) and (\ref{eq:7}) $\widehat{f}$ in  (\ref{eq:5}) is the L.t. of
\begin{eqnarray}
\label{eq:10}
f(x_1,\ldots,x_p) = (2\pi)^{-p} \int_{\cC_p} B(y^{-1}_{1}, \ldots,y^{-1}_{p}) \prod^{p}_{j=1} g_j(x_j,y_j) d\varphi_j
\end{eqnarray}
with $y_j = re^{i\varphi_j}$, \ $-\pi < \varphi \leq \pi$, \ $r^{-1}_{B} < r < 1$, \ $g_j$ from (\ref{eq:8}). 
\end{theorem}

\begin{proof}
The integral in (\ref{eq:10}) is evaluated by
\begin{eqnarray*}
\begin{array}{l}
(2\pi)^{-p}  \displaystyle\int_{\cC_p} \left(\sum b(m_1,\ldots,m_p) \prod^{p}_{j=1} y^{-m_j}_{j} \right)\left( \sum \prod^{p}_{j=1} g_{jn_j}(x_j)y_{j}^{n_j}\right)d\varphi_1\ldots d\varphi_p\\
= \sum b(n_1,\ldots,n_p) \prod^{p}_{j=1} g_{jn_j}(x_j)
\end{array}
\end{eqnarray*}
and this series has the L.t. from (\ref{eq:5}). 
\end{proof}

Some further remarks: With
\begin{eqnarray}
\label{eq:11}
G_j(x_j,y_j) = \int^{x_j}_{0} g_j(\xi, y_j)d\xi
\end{eqnarray}
instead of the $g_j$ in (\ref{eq:10}), the corresponding representation arises for
\begin{eqnarray}
\label{eq:12}
F(x_1,\ldots,x_p) = \int^{x_1}_{0}\ldots \int^{x_p}_{0} f(\xi_1,\ldots,\xi_p) d\xi_1 \ldots d\xi_p .
\end{eqnarray}
If the series in (\ref{eq:8}) are absolutely convergent for all $y \in \bC$ then additionally 
\begin{eqnarray}
\label{eq:13}
\lim_{t\to \infty} h_j(t) = 0
\end{eqnarray}
is supposed to hold. Then the rhs of (\ref{eq:9}) is the L.t. of $g_j(x,y)$ for any fixed $y$ and all sufficiently large $t$.

In some cases the functions $g_{j_0}$ and their L.t. $\widehat{g}_{j_0}$ are known from univariate marginal distributions apart from some scale factors. If the functions $u_j = h_j(t)$ are explicitly invertible then
\begin{eqnarray}
\label{eq:14}
B(u_1,\ldots,u_p) = \frac{\widehat{f}\left(h^{-1}_{1}(u_1),\ldots,h^{-1}_{p}(u_p)\right)}{\prod^{p}_{j=1}\widehat{g}_{j_0}\left(h^{-1}_{j}(u_j)\right)}
\end{eqnarray}
can sometimes be found easily from the given $\widehat{f}$.

\section{\large Three representations for the $\Gamma_p(\alpha,\Sigma,\Delta)$--Laplace transform \newline
and the related generators}
\label{section:3}
With any $v > 0$ we define
\begin{eqnarray}
\label{eq:15}
z_j= (1+v^{-1}t_j)^{-1}, \ t_j \geq 0, \ u_j = 1-z_j = v^{-1}t_jz_j, \ \omega_j = z_j - u_j, &&\nonumber\\[-2ex]
&&\\[-2ex]
Z = \diag(z_1,\ldots, z_p), \ U = \diag(u_1, \ldots, u_p), \ \Omega = \diag(\omega_1, \ldots,\omega_p).&&\nonumber
\end{eqnarray}
The scale factor $v$ is introduced to obtain $\|B\| < 1$ for the matrices $B$ defined in (\ref{eq:20}) below and to effect the convergence of some series expansions. For a more general scaling see remarks following theorem~2 in section~4.

From the relations 
\begin{eqnarray}
\label{eq:16}
v^{-1}T = UZ^{-1}, \ I_p = Z + U, \ \Omega = Z-U,
\end{eqnarray}
it follows for the matrices $I + \Sigma T$ in the L.t. (\ref{eq:1}):
\begin{eqnarray}
\label{eq:17}
I + \Sigma T = I + v\Sigma U Z^{-1}= (Z + v\Sigma U)Z^{-1}
\end{eqnarray}
and

\newcounter{appeqn}
\setcounter{appeqn}{1}

\begin{eqnarray*}
\label{eq:18}
\hspace*{2cm} Z + v\Sigma U = \left\{ \begin{array}{l@{\hspace{2.2cm}}r}
I + (v\Sigma - I)U, & (\arabic{equation}\alph{appeqn})\\
\stepcounter{appeqn}
v\Sigma (I + (v^{-1} \Sigma^{-1} - I)Z), & (\arabic{equation}\alph{appeqn})\\
\stepcounter{appeqn}
\frac 12 (I + v\Sigma)( I + (2(I + v\Sigma)^{-1} - I)\Omega),\qquad  & (\arabic{equation}\alph{appeqn})
\end{array}\right.
\end{eqnarray*}

\setcounter{appeqn}{1}
\stepcounter{equation}

and therefore 
\begin{eqnarray}
\label{eq:19}
|I + \Sigma T|^{-\alpha} = c^\alpha |Z|^\alpha |I + BY|^{-\alpha}
\end{eqnarray}
with
\begin{eqnarray*}
\label{eq:20}
\hspace*{2.5cm}
\begin{array}{lll@{\hspace{2.6cm}}r}
Y = U, \quad & B = v\Sigma - I, & c = 1, & (\arabic{equation}\alph{appeqn})\\
\stepcounter{appeqn}
Y = Z, & B = (v\Sigma)^{-1} - I, & c = |I+B|, & (\arabic{equation}\alph{appeqn})\\
\stepcounter{appeqn}
Y = \Omega, & B = 2(I +v\Sigma)^{-1}-I, \quad & c = |I+B|.\quad & (\arabic{equation}\alph{appeqn})
\end{array}
\end{eqnarray*}

\setcounter{appeqn}{1}
\stepcounter{equation}

It should be noticed that $\|B\|<1$ in (\ref{eq:20}c) for every $v >0$ and $\Sigma >0$.

Now, using (\ref{eq:16}), by a straightforward calculation the L.t. in (\ref{eq:1}) can be represented by 
\begin{eqnarray*}
\label{eq:21}
\lefteqn{\widehat{f}(t_1,\ldots,t_p; \alpha,\Sigma,\Delta) =\nonumber}\\[3ex]
\hspace*{1cm}&&
\left\{
\begin{array}{r@{\hspace{1.7cm}}r}
|Z|^\alpha|I+BU|^{-\alpha} \mathrm{etr}(-(I+B)U(I+BU)^{-1}\Delta), &(\arabic{equation}\alph{appeqn})\\
\stepcounter{appeqn}
 |I+B|^\alpha \mathrm{etr}(-\Delta)|Z|^\alpha |I+BZ|^{-\alpha}\mathrm{etr}(Z(I+BZ)^{-1}(I+B)\Delta), &(\arabic{equation}\alph{appeqn})\\
\stepcounter{appeqn}
 |I+B|^\alpha \mathrm{etr}(-\frac 12 \Delta (I-B))& (\arabic{equation}\alph{appeqn})\\
\cdot |Z|^\alpha |I+B\Omega|^{-\alpha}\mathrm{etr}(\frac 12 \Omega(I+B\Omega)^{-1}(I+B)\Delta(I-B)),  & 
\end{array} 
\right.
\end{eqnarray*}

\setcounter{appeqn}{1}
\stepcounter{equation}

with the corresponding matrices $B$ from (\ref{eq:20}) and $Z,U,\Omega$ from (\ref{eq:15}).

For the former series expansions the following relations were used:

Laplace  transform  $\widehat{f}(t)$: \hspace{1.2cm} $f(x)$: \hspace{1.5cm} \mbox{$F(x)\,=\,\int_0^x f(\xi)d\xi$:}\\ 
\rule[.3mm]{14cm}{.3mm}
\begin{eqnarray*}
\label{eq:22}
\hspace*{1cm}
\begin{array}{l@{\hspace{3.5cm}}l@{\hspace{2cm}}l@{\hspace{2.2cm}}r}
z^\alpha u^n        &  vg_{\alpha+n}^{(n)}(vx)       & G_{\alpha+n}^{(n)}(vx) &(\arabic{equation}\alph{appeqn})\\
\stepcounter{appeqn}
z^{\alpha+n}        & vg_{\alpha+n}(vx)             & G_{\alpha+n}(vx)& (\arabic{equation}\alph{appeqn})\\
\stepcounter{appeqn}
z^{\alpha}\omega^n &  vh_{\alpha,n}(vx)             & H_{\alpha,n}(vx)  & (\arabic{equation}\alph{appeqn})
\end{array}
\end{eqnarray*}

\setcounter{appeqn}{1}
\stepcounter{equation}

where $z = (1+v^{-1}t)^{-1}$, $g_{\alpha+n}(x) = e^{-x} x^{\alpha-1+n}/ \Gamma (\alpha +n)$, \newline 
\[
g_{\alpha+n}^{(n)}(x) = \frac{d^n}{dx^n} \  g_{\alpha+n}(x) = {\alpha - 1 + n \choose n}^{-1} L^{(\alpha-1)}_{n}(x) g_\alpha(x)
\]
  with the generalized Laguerre polynomials $L^{(\alpha-1)}_{n}$ and 
\[
  h_{\alpha,n}(x) = (-1)^n {\alpha - 1 + n \choose n}^{-1} L^{(\alpha-1)}_{n}(2x) g_\alpha(x).
\]
The last identity is verified by L.t.

The following bounds are derived from (22.14.13) in Abramowitz and Stegun (1965):
\begin{eqnarray}
\label{eq:23}
\left| g_{\alpha+n}^{(n)} (x) \right| \leq 
\left\{
\begin{array}{l@{\hspace{.5cm}}l}
e^{x/2} g_\alpha(x), & \alpha \geq 1\\
2n\alpha^{-1} e^{x/2} g_\alpha (x), & 0 < \alpha < 1
\end{array}
\right\}
\end{eqnarray}
\begin{eqnarray}
\label{eq:24}
\left| h_{\alpha,n}(x)\right| \leq 
\left\{
\begin{array}{l@{\hspace{.5cm}}l}
x^{\alpha-1}/\Gamma(\alpha), & \alpha \geq 1\\
2nx^{\alpha-1}/\Gamma(\alpha+1), & 0 < \alpha < 1
\end{array}
\right\},
\end{eqnarray}
matching with the conditions in (\ref{eq:7}).

The following generators (generating functions) with the $\Gamma(\alpha+n)$--cdf $G_{\alpha + n}(x)$ are required for the formulas in \mbox{theorem~2:}
\begin{eqnarray*}
\label{eq:25}
\hspace*{1cm}
F_{\alpha}(x,y)=
\left\{
\begin{array}{l@{\hspace{.5cm}}l@{\hspace{2cm}}r}
\sum^{\infty}_{n=0} G^{(n)}_{\alpha+n}(x)y^n = \frac 1{1-y}\  G_\alpha \left(x,\frac y{y-1}\right), & |y|<1, & (\arabic{equation}\alph{appeqn})\\
\stepcounter{appeqn}
\sum^{\infty}_{n=0} G_{\alpha+n}(x)y^n =  G_\alpha (x,y), & y \in \bC, & (\arabic{equation}\alph{appeqn})\\
\stepcounter{appeqn}
\sum^{\infty}_{n=0} H_{\alpha,n}(x)y^n = \frac 1{1+y} \ G_\alpha \left(x,\frac {2y}{y+1}\right), & |y|<1 & (\arabic{equation}\alph{appeqn})
\end{array}
\right. 
\end{eqnarray*}

\setcounter{appeqn}{1}
\stepcounter{equation}

The identities (a) and (c) are verified by the L.t. of $f_\alpha (x,y) = \frac\partial{\partial x} F_\alpha(x,y)$. A short calculation shows
\begin{eqnarray}
\label{eq:26}
G_{\alpha}(x,y)=
\left\{
\begin{array}{l@{\hspace{.5cm}}l@{\hspace{.5cm}}ll}
\frac 1{1-y} \left( G_\alpha(x) - y^{1-\alpha}e^{(y-1)x}\ G_\alpha(xy)\right), & y\neq 1, & \alpha > 0\\
\frac 1{1-y} \left( G_{\alpha-1}(x) - y^{1-\alpha}e^{(y-1)x}\ G_{\alpha-1}(xy)\right), & \alpha \geq 1, & G_0 := 1\\
xg_\alpha(x) + (1+x-\alpha) \ G_\alpha(x), & y=1 &
\end{array}
\right.
\end{eqnarray}
and 
\begin{eqnarray*}
\label{eq:26*}
g_\alpha (x,y) = \frac\partial{\partial x}\ G_{\alpha}(x,y)=
\left\{
\begin{array}{l@{\hspace{.5cm}}l@{\hspace{.5cm}}l}
g_\alpha(x) + y^{1-\alpha} e^{(y-1)x} \ G_\alpha(xy), & \alpha > 0\\
y^{1-\alpha} e^{(y-1)x} \ G_{\alpha-1}(xy), & \alpha \geq 1
\end{array}
\right\}.
\end{eqnarray*}
The functions $F_\alpha(x,y)$ are especially simple for $\alpha \in \bN$ since $G_\alpha (z) = 1 - e^{-z} \sum^{\alpha-1}_{j=0}{z^j}/{j!}$, $\alpha \in \bN$.

Besides, 
\[
G_{k+1/2}(z) = \mathrm{erf}(z^{1/2}) - e^{-z}\sum^{k}_{j=1} \frac{z^{j-1/2}}{\Gamma(j + 1/2)}, \ k\in \bN_0.
\]
The following simple lemma is used for the proof of theorem~2.

\begin{lemma}
\label{lemma:1}
If $B$ is a symmetrical $p\times p$--matrix with $\|B\|< 1$ and \\
$Y = \diag (y_1,\ldots,y_p)$ then the power series expansion
\[
|I + BY|^{-\alpha} = \sum b(n_1,\ldots, n_p) \prod^{p}_{j=1} y^{n_j}_{j}
\]
is absolutely convergent for $\max |y_j| < r_B = \|B\|^{-1}$.
\end{lemma}

This follows from $\sum_{(n)} |b(n_1,\ldots,n_p)| = O(\vartheta^n)$ with any $\vartheta > \|B\|$, which has been already shown in (2.1.16) $\ldots$ (2.1.18) in Royen (1991) (with the notation  $-C$ instead of $B$).

\section{\large The integral representations}
\label{section:4}
In theorem~2 below the functions $F(x_1,\ldots,x_p; \alpha_1,\ldots,\alpha_n, \Sigma_1,\ldots,\Sigma_n, \Delta_1,\ldots,\Delta_n)$ from (\ref{eq:2}) are represented by three different integrals over $\cC_p= (-\pi,\pi]^p$. Together with the generators $F_\alpha$ from (\ref{eq:25}), $\alpha = \sum^{n}_{k=1} \alpha _k$, the following matrices are used with a scale factor $v$ to enforce $\|B_k\|< 1$:
\begin{eqnarray*}
\label{eq:27}
\hspace*{.5cm}
\begin{array}{l@{\hspace{.3cm}}l@{\hspace{.3cm}}l@{\hspace{1.2cm}}r}
B_k = v\Sigma_k - I, & D_k = \Delta _k(I+B_k), &F_\alpha \ \mathrm{from \ (\ref{eq:25}a)}, & (\arabic{equation}\alph{appeqn})\\
\stepcounter{appeqn}
B_k = (v\Sigma_k)^{-1} - I, & D_k =(I+B_k)\Delta_k, & F_\alpha \ \mathrm{from \ (\ref{eq:25}b)}, & (\arabic{equation}\alph{appeqn})\\
\stepcounter{appeqn}
B_k = 2(I+v\Sigma_k)^{-1} - I, & D_k = \frac 12 (I+B_k)\Delta_k(I-B_k), & F_\alpha \ \mathrm{from \ (\ref{eq:25}c)}. & (\arabic{equation}\alph{appeqn})
\end{array}
\end{eqnarray*}

\setcounter{appeqn}{1}
\stepcounter{equation}

Furthermore, we define $\lambda_{\max} = \max \|\Sigma_k\|$ , \ $\lambda_{\min}^{-1} = \max\|\Sigma_k^{-1}\|$, \ $y_j =re^{i\varphi_j}$, \\
$-\pi < \varphi_j \leq \pi$, \ $Y = \diag (y_1,\ldots,y_p)$, 
\[
K = K(y_1,\ldots,y_p) = \prod^{n}_{k=1}\ \mathrm{etr}(\pm(Y+B_k)^{-1}D_k) |I+B_k Y^{-1}|^{-\alpha_k},
\]
where the negative sign occurs only with $B_k, D_k $ from  (\ref{eq:27}a), and \\
$\cF_\alpha d\varphi = \prod^{p}_{j=1}F_\alpha(vx_j, y_j)d\varphi_j$.

\begin{theorem}
\label{theorem:2}
With the above notations the functions $F$ from (\ref{eq:2}) are respresentable by each of the following three integrals:
\begin{eqnarray}
\label{eq:28}
(2\pi)^{-p} \int_{\cC_p} K \cF_\alpha d\varphi,
\end{eqnarray}
$F_\alpha$ from (\ref{eq:25}a), $B_k,D_k$ from (\ref{eq:27}a), $\|B_k\|<1$ if $v < 2\lambda^{-1}_{\max}$, $ \max \|B_k\|< r < 1$,
\begin{eqnarray}
\label{eq:29}
\left( \prod^{n}_{k=1}\ \mathrm{etr} (-\Delta_k)|I+B_k|^{\alpha_k}\right) (2\pi)^{-p} \int_{\cC_p}K\cF_\alpha d\varphi,
\end{eqnarray}
$F_\alpha$ from (\ref{eq:25}b), $B_k,D_k$ from (\ref{eq:27}b), $\|B_k\|<1$ if $v > \frac 12\lambda^{-1}_{\min}$, $ \max \|B_k\|< r$,
\begin{eqnarray}
\label{eq:30}
\left( \prod^{n}_{k=1}\ \mathrm{etr} \left(-\textstyle \frac 12 \Delta_k(I-B_k)\right) |I+B_k|^{\alpha_k}\right) (2\pi)^{-p} \int_{\cC_p}K\cF_\alpha d\varphi,
\end{eqnarray}
$F_\alpha$ from (\ref{eq:25}c), $B_k,D_k$ from (\ref{eq:27}c), $v >0$, $ \max \|B_k\|< r < 1$.
\end{theorem}

\begin{proof}
Because of lemma~1 the assumptions of theorem~1 are satisfied with \mbox{$\widehat{g}_{j_0}(t) = z^\alpha_j = (1+v^{-1}t_j)^{-\alpha}$} and $h_j(t_j)$ corresponding to $z_j$ or $u_j = v^{-1}t_jz_j = 1-z_j$ or \mbox{$\omega_j = z_j-u_j$} respectively. The functions $\widehat{g}_{j_0}(t)(h_j(t))^n$ are the L.t. of the functions in the second column of (\ref{eq:22}) from which type (a) and (c) have the bounds in (\ref{eq:23}),  (\ref{eq:24}), satisfying the condition (\ref{eq:7}) for theorem~1. The series $\sum^{\infty}_{n=0} G_{\alpha+n} (x) y^n = G_\alpha(x,y)$ in (\ref{eq:25}b) is absolutely convergent for every $y \in \bC$. Thus, all $r > \max \|B_k\|$ are admissible in (\ref{eq:29}). In (\ref{eq:30}) we have $\max \|B_k\|< 1$ for every $v >0$. Hence, theorem~1 together with the respresentations of the L.t. in (\ref{eq:21}) implies (\ref{eq:28}), (\ref{eq:29}) and (\ref{eq:30}).
\end{proof}

The univariate case of (\ref{eq:29}) provides 
\begin{eqnarray}
\label{eq:31}
\lefteqn{F(x;\alpha_1,\ldots,\alpha_n, \sigma^{2}_{1},\ldots,\sigma^{2}_{n}, \delta^2_1, \ldots, \delta^2_n)=}\nonumber\\[-2ex]
&& \\[-2ex]
&&\left( v^{-\alpha} \prod_{k=1}^n \sigma_k^{-2\alpha_k}e^{-\delta^2_k}\right)
 \frac 1\pi \int^{\pi}_{0} \cR e 
\left\{ 
G_\alpha ( vx, e^{i\varphi})\prod^n_{k=1}
\frac{\exp \left( \delta^2_k / (1 +  v \sigma^2_k (e^{i\varphi}-1))\right)}
{\left( 1+ (v^{-1} \sigma_k^{-2}-1)  e^{-i\varphi} \right)^{\alpha_k}}
\right\} d\varphi\nonumber
\end{eqnarray}
with $2v > \max \sigma^{-2}_k$, $r=1$, $G_\alpha$ from (\ref{eq:26}). With $p=1$ similar formulas arise from (\ref{eq:28}) or (\ref{eq:30}).

The cdf of a quadratic form $\frac 12 x'Ax$ with $T'AT = \diag (\lambda_1,\ldots,\lambda_n) \geq 0$ of rank $q$ and a $\cN(\mu, \sigma^2 I_n)$--random vector $x$ is a special case of (\ref{eq:31}) with $\alpha_k=\frac 12$, $\sigma^2_k = \lambda_k \sigma^2$ and non--centrality parameters $\delta^2_k = \frac 12 \mu^{*2}_{k}/\sigma^2$, $k = 1,\ldots,q$, $\mu^*=T'\mu$.

Some further remarks: In (\ref{eq:29}) also $\|B_k\|>1 $ is allowed since every $r = \|Y \|> \max \|B_k\|$ is admissible, which entails $\max \|B_k Y^{-1}\| < 1$.

With $\vartheta = \lambda_{\max}/\lambda_{\min}$ it follows with special values of $v$: 
\begin{eqnarray*}
\max \|B_k\| \leq \frac{\vartheta-1}{\vartheta+1} \ &\mbox{ in (\ref{eq:28}) with } \ & v = 2(\lambda_{\min} + \lambda_{\max})^{-1}, \\
\max \|B_k\| \leq \frac{\vartheta-1}{\vartheta+1} \ &\mbox{ in (\ref{eq:29}) with }  \ & v = \frac 12(\lambda^{-1}_{\min} + \lambda^{-1}_{\max}), 
\end{eqnarray*}
but
\begin{eqnarray*}
\max \|B_k\| \leq \frac{\sqrt{\vartheta}-1}{\sqrt{\vartheta}+1} \ \mbox{ in (\ref{eq:30}) with } \ v = (\lambda_{\min} \lambda_{\max})^{-1/2}.
\end{eqnarray*}

More generally, the scale factor $v = w^2$ can be replaced by a scale matrix \mbox{$W^2 = \diag (w^{2}_{1},\ldots,w^{2}_{p}) > 0$.} Then with $T_w = W^{-1}TW^{-1}$, $\Sigma_w = W\Sigma W$, $\Delta_w = W\Delta W^{-1}$ the L.t. (\ref{eq:1}) equals
\begin{eqnarray}
\label{eq:32}
|I + \Sigma_w T_w|^{-\alpha} \ \mathrm{etr}(-\Sigma_wT_w(I+\Sigma_wT_w)^{-1} \Delta_w).
\end{eqnarray}
Consequently, besides the substitutions $v\Sigma_k \to W\Sigma_kW$, $\Delta_k \to W\Delta_kW^{-1}$, the matrices $I + B_k$ in theorem~2 must be replaced by $W \Sigma _kW$, $(W\Sigma_kW)^{-1}$ and $2(I + W\Sigma_kW)^{-1}$ respectively, and the generators $F_\alpha(vx_j, y_j)$ by $F_\alpha(w_j^2 x_j, y_j)$.

In particular for a single $\Gamma_p(\alpha, \Sigma,\Delta)$--distribution this more general scaling can be used to minimize $\|B\|$ or for a "natural scaling" i.e. to standardize $I +B$ to a correlation matrix. However, $\|B\|<1$ must be taken into account in (\ref{eq:28}), whereas this condition is satisfied in (\ref{eq:30}) for every scaling. It was shown in Royen (1991) that natural scaling can always be accomplished also in $I+B = 2(I + W \Sigma W)^{-1}$ by a unique $W^2$.

\section{\large Representations of the $\Gamma_p(\alpha,\Sigma,\Delta)$ distribution function 
by $(p-1)$--variate integrals}
\label{section:5}
For a single $\Gamma_p(\alpha,\Sigma,\Delta)$--cdf it is always possible to perform the integration over a single variable $\varphi_j$ within the integrals from theorem~2.

We use the following functions 
\begin{eqnarray}
\label{eq:33}
\begin{array}{l}
\displaystyle
\cG_\alpha (x,y) = e^{-y} \sum^{\infty}_{n=0} G_{\alpha +n}(x) \frac{y^n}{n!} = \sum^{\infty}_{n=0} G^{(n)}_{\alpha +n}(x) \frac{(-y)^n}{n!},\\
x,y \in \bC, \ G_{\alpha+n}, \ G^{(n)}_{\alpha+n} \ \mbox{ from (\ref{eq:22}), \ and }\\
\cG^*_\alpha(x,y) = e^y \cG_\alpha (x,y). 
\end{array}
\end{eqnarray}
For positive half integers $\alpha = 1/2 + k$ these functions can also be computed by the erf--function and a sum of $k$ terms which are essentially given by the modified Bessel functions $I_{j-1/2}(2(xy)^{1/2})$, $j=1,\ldots,k$, (see e.g. Royen (1995) or (2006)).

Now let be $W^2 = \diag(w^2_1,\ldots,w^2_p)$ a general scale matrix, \\
$Y = \diag (y_1,\ldots,y_p)$, $y_j = re^{i\varphi_j}$, $-\pi < \varphi_j \leq \pi$,
\begin{eqnarray*}
\label{eq:34}
\hspace*{1.7cm}
B=
\begin{pmatrix}
B_{pp} & b_p\\ b'_p & b_{pp} \end{pmatrix} = \left\{
\begin{array}{l@{\hspace{4cm}}r}
W\Sigma W - I, & (\arabic{equation}\alph{appeqn})\\
\stepcounter{appeqn}
(W \Sigma W)^{-1} - I , & (\arabic{equation}\alph{appeqn})\\
\stepcounter{appeqn}
2(I + W \Sigma W)^{-1} - I, \quad & (\arabic{equation}\alph{appeqn})
\end{array}
\right.
\end{eqnarray*}

\setcounter{appeqn}{1}
\stepcounter{equation}

\begin{eqnarray*}
\label{eq:35}
\hspace*{.7cm}
D=
\begin{pmatrix}
D_{pp} & d_p\\ d^p& d_{pp} 
\end{pmatrix} = 
\left\{
\begin{array}{l@{\hspace{1cm}}r}
W\Delta\Sigma W , &(\arabic{equation}\alph{appeqn})\\
\stepcounter{appeqn}
W^{-1} \Sigma^{-1}\Delta W^{-1}, & (\arabic{equation}\alph{appeqn})\\
\stepcounter{appeqn}
2(I + W \Sigma W)^{-1} W\Delta W^{-1}(I-(I+W \Sigma W)^{-1}), \quad & (\arabic{equation}\alph{appeqn})
\end{array}\right.
\end{eqnarray*}

\setcounter{appeqn}{1}
\stepcounter{equation}

\begin{equation}
\label{eq:36}
y_0 = y_0(y_1,\ldots,y_{p-1}) = b'_p(Y_{pp}+B_{pp})^{-1}b_p-b_{pp}
\end{equation}
\begin{equation}
\label{eq:37}
q = q(y_1,\ldots, y_{p-1}) = (b'_p(Y_{pp}+B_{pp})^{-1}, -1) D 
 {(Y_{pp}+B_{pp})^{-1} b_p\choose -1}
\end{equation}
and 
\[
K_\alpha = K_\alpha(y_1,\ldots, y_{p-1}) = \mathrm{etr}(\pm(Y_{pp}+B_{pp})^{-1}D_{pp})|I + B_{pp}Y^{-1}_{pp}|^{-\alpha},
\]
where the negative sign is only taken for $B_{pp}$ from (\ref{eq:34}a).

\begin{theorem}
\label{theorem:3} 
With the above notations the $\Gamma_p(\alpha,\Sigma,\Delta)$--cdf  $F(x_1,\ldots,x_p;\alpha, \Sigma, \Delta)$ is given by each of the following three integrals:
\begin{equation}
\label{eq:38}
\frac{1}{(2\pi)^{p-1}} \int_{\cC_{p-1}} \cG_\alpha \left( \frac{w_p^2 x_p}{1-y_0}, \frac q{1-y_0}\right) K_\alpha \prod^{p-1}_{j=1} \frac 1{1-y_j} \  G_\alpha \left( w^2_jx_j, \frac{y_j}{y_j-1}\right) d\varphi_j,
\end{equation}
$B$ from (\ref{eq:34}a), $D$ from (\ref{eq:35}a), $\|B\|<r<1$,
\begin{eqnarray}
\label{eq:39}
\frac{\mathrm{etr}(-W \Delta W^{-1})}{|W\Sigma W|^\alpha}& \cdot & \frac 1{(2\pi)^{p-1}} \int_{\cC_{p-1}} (1-y_0)^{-\alpha} \ \cG^*_\alpha
\left( (1-y_0) w^2_p x_p, \frac q{1-y_0}\right)\nonumber\\[-2ex]
&&\\[-2ex]
& \cdot & K_\alpha \prod^{p-1}_{j=1} G_\alpha (w^2_j x_j,y_j)d\varphi_j,\nonumber
\end{eqnarray}
$B$ from (\ref{eq:34}b), $D$ from (\ref{eq:35}b), $\|B\|<r$,
\begin{eqnarray}
\label{eq:40}
\lefteqn{\frac{2^{\alpha p}\mathrm{etr}(-\frac 12 W \Delta W^{-1}(I-B))}{|I + W\Sigma W|^\alpha} \cdot  \frac 1{(2\pi)^{p-1}} \cdot \int_{\cC_{p-1}} \exp \left( \frac q {1-y_0}\right) .\nonumber}\\[-2ex]
&&\\[-2ex]
&&(1-y_0)^{-\alpha} \cG_\alpha \left( \frac{1-y_0}{1+y_0}  w^2_p x_p, \frac {2q} {1-y^2_0}\right)
 K_\alpha \prod^{p-1}_{j=1} \frac{ 1} { 1+y_j} \ G_\alpha \left(w^2_j x_j,\frac{2y_j}{y_j+1}\right)d\varphi_j,\nonumber
\end{eqnarray}
$B$ from (\ref{eq:34}c), $D$ from (\ref{eq:35}c), $\|B\|<r<1$.
\end{theorem} 
For the proof of theorem~3 the following two lemmas are required.

\begin{lemma}
\label{lemma:2}
With $Y = \diag(y_1,\ldots,y_p)$, $y_j = re^{i\varphi_j}$, $\|B\| < r$, $B,D, y_0, q$ from (\ref{eq:34}), (\ref{eq:35}), (\ref{eq:36}), (\ref{eq:37}) the following decomposition is obtained
\begin{eqnarray}
\label{eq:41}
\lefteqn{\mathrm{etr} ((Y+B)^{-1}D) |Y+B|^{-\alpha}\nonumber}\\[-2ex]
&&\\[-2ex]
&=&\mathrm{etr} \left((Y_{pp}+ B_{pp})^{-1} D_{pp}\right) |Y_{pp} + B_{pp}|^{-\alpha} \exp\left(\frac q{y_p-y_0}\right) (y_p-y_0)^{-\alpha}.\nonumber
\end{eqnarray}
\end{lemma}

\begin{proof}
From frequently used formulas for $p\times p$--matrices, (see e.g. complements and problems 2.4, 2.7 in chapter 1b of Rao (1973)) it follows for \\
$A = Y+B= \scriptstyle{\begin{pmatrix}  A_{pp} & b_{p} \\ b'_p & \qquad y_p+b_{pp}\end{pmatrix}}$:
\begin{eqnarray*}
|A| &=& |A_{pp}| (y_p + b_{pp}- b'_p A^{-1}_{pp} b_p) = |Y_{pp} + B_{pp}| (y_p-y_0),\\
A^{-1}&=& 
\begin{pmatrix}
A^{-1}_{pp} + \frac 1{y_p-y_0} A^{-1}_{pp} b_p b'_p A^{-1}_{pp} \qquad & -\frac 1{y_p-y_0} A^{-1}_{pp} b_p\\
 -\frac1{y_p-y_0} b'_p A^{-1}_{pp} & \frac 1{y_p- y_0} 
\end{pmatrix}
\end{eqnarray*}
and 
\begin{eqnarray*}
\lefteqn{\mathrm{trace}(A^{-1}D)}\\
&=& \mathrm{trace} \left(A^{-1}_{pp} D_{pp}+ \frac 1{y_p-y_0} \left(A^{-1}_{pp} b_p b'_p A^{-1}_{pp} D_{pp} - A^{-1}_{pp} b_p d^p \right)\right) +\frac 1{y_p-y_0}(d_{pp}-b'_p A_{pp}^{-1}d_p)\\
&=& \mathrm{trace} (A^{-1}_{pp} D_{pp}) + \frac q{y_p-y_0}, \ \mbox{which implies (\ref{eq:41})}.
\end{eqnarray*}
\end{proof}

\begin{lemma}
\label{lemma:3}
Let be $q$ any number, $S_r = \{y \in \bC \big| |y| = r\}$, $y_0$ any number with $|y_0| < r$, then with $F_\alpha$ from (\ref{eq:25}), $\cG_\alpha, \cG^*_\alpha$ from (\ref{eq:33}), and the negative sign in $\pm q$ only for (\ref{eq:42}a)
\begin{eqnarray*}
\label{eq:42}
\lefteqn{\frac 1{2\pi i} \oint_{S_r}\mathrm{etr}\left( \frac {\pm q}{y-y_0}\right) F_\alpha (x,y) (y-y_0)^{-\alpha} y^{\alpha -1} dy\nonumber }\\[2ex]
&&= \left\{
\begin{array}{ll@{\hspace{1.2cm}}r}
\cG_\alpha \left( \frac x{1-y_o}, \frac q{1-y_0}\right) & F_\alpha \ \mbox{from (\ref{eq:25}a), $r < 1$}, \quad  & (\arabic{equation}\alph{appeqn})\\
\stepcounter{appeqn}
(1-y_0)^{-\alpha} \ \cG^*_\alpha \left( (1-y_0)x, \frac q{1-y_0}\right), & F_\alpha \ \mbox{from (\ref{eq:25}b)}, & (\arabic{equation}\alph{appeqn})\\
\stepcounter{appeqn}
 \exp\left(\frac q{1-y_0}\right) (1-y_0)^{-\alpha} \ \cG_\alpha \left( \frac{1-y_0}{1+y_0} x, \frac{2q}{1-y^2_0} \right), \quad    & F_\alpha \ \mbox{from (\ref{eq:25}c), $r < 1$}, &(\arabic{equation}\alph{appeqn})
\end{array}
\right. 
\end{eqnarray*}

\setcounter{appeqn}{1}
\stepcounter{equation}

\end{lemma}

\begin{proof}
It is sufficient to verify (\ref{eq:42}) for the corresponding derivatives $f_\alpha = \frac \partial{\partial x}F_\alpha$. At first, (\ref{eq:42}a) is shown: 

With $F_\alpha$ from (\ref{eq:25}a) and the binomial series for $(1-y_0/y)^{-(\alpha+n)}$ we obtain
\begin{eqnarray*}
\lefteqn{\frac 1{2 \pi i} \oint_{S_r} f_\alpha (x,y) (y-y_0)^{-(\alpha+n)}y^{\alpha-1} dy}\\
&=& \frac 1{2\pi i} \oint_{S_r} \left(\sum^{\infty}_{m=0} g^{(m)}_{\alpha+m}(x) y^m\right)\left( \sum^{\infty}_{k=0}{\alpha+n+k-1\choose k}\left(\frac {y_0}y\right)^k\right) y^{-n-1}dy.
\end{eqnarray*}
With $z =(1+t)^{-1}$, $u = tz$, the last integral has the L.t.
\begin{eqnarray*}
\lefteqn{\frac 1{2 \pi i} \oint_{S_r} z^\alpha \left( \sum^{\infty}_{m=0} (uy)^m\right) \cdot \left( \sum^{\infty}_{k=0}{\alpha+n+k-1\choose k}\left(\frac {y_0}y\right)^k\right) y^{-n-1}dy}\\
&=& z^\alpha \sum_{m=n+k} u^m {\alpha+n+k-1\choose k} y_0^k = z^\alpha u^n (1-uy_0)^{-(\alpha+n)}.
\end{eqnarray*}
Multiplication by $(-q)^n/n!$ and summation over $n$ leads to the L.t.
\begin{eqnarray*}
\frac {z^\alpha}{(1-uy_0)^\alpha} \exp\left( - \frac{qu}{1-uy_0}\right) = \frac 1{(1+(1-y_0)t)^\alpha} \exp\left( - \frac{\frac q{1-y_0}(1-y_0)t}{1+(1-y_0)t}\right)
\end{eqnarray*}
and this is the L.t. of $\frac \partial{\partial x} \ \cG_\alpha \left( \frac x{1-y_0}, \frac q{1-y_0}\right)$.

To verify (\ref{eq:42}b) we obtain with $F_\alpha$ from (\ref{eq:25}b):

\begin{eqnarray*}
\lefteqn{\frac 1{2 \pi i} \oint_{S_r} f_\alpha (x,y)(y-y_0)^{-(\alpha+n)} y^{\alpha-1}dy}\\
&&=\frac 1{2\pi i} \oint_{S_r} \left(\sum^{\infty}_{m=0} g_{\alpha+m}(x) y^m\right)\left( \sum^{\infty}_{k=0}\frac{\Gamma(\alpha+n+k)}{\Gamma(\alpha+n)k!}\left(\frac {y_0}y\right)^k\right) y^{-n-1}dy\\
&&=\sum_{m=n+k}g_{\alpha+m}(x)\frac{\Gamma(\alpha+n+k)}{\Gamma(\alpha+n)k!} y_0^k \ = \ g_{\alpha+n}(x) e^{xy_0}\\
&& =(1-y_0)^{-(\alpha+n)} (1-y_0) g_{\alpha+n}((1-y_0)x). 
\end{eqnarray*}
Multiplication by $q^n/n!$ and summation provides $(1-y_0)^{-\alpha} \frac \partial {\partial x} \ \cG^* \left((1-y_0)x , \frac q{1-y_0}\right)$.

(\ref{eq:42}c) can be shown by L.t. in a similar way as (\ref{eq:42}a).

\paragraph{\sc Proof of theorem~3. } Without loss of generality $y_p$ is selected from the variables $y_j = re^{i\varphi_j}$ in $Y = diag(y_1,\ldots,y_p)$ with any fixed $r > \|B\|$. If $y_p$ is replaced by a variable $y$ with any $|y|$ then the equation
\begin{eqnarray*}
|Y+B| = |Y_{pp} + B_{pp}|(y_p+b_{pp} -b'_p(Y_{pp} + B_{pp})^{-1} b_p) = 0
\end{eqnarray*}
has always a unique solution
\begin{eqnarray*}
y = y_0 = b'_p(Y_{pp} + B_{pp})^{-1} b_p-b_{pp}
\end{eqnarray*}
with $|y_0| < r$ since $\|B_{pp}\| \leq \|B\|$.

Hence, with lemma~2 and lemma~3, theorem~3 is obtained by integration over $\varphi_p$ in the integrals of theorem~2 with $n=1$.

\end{proof}


%


%
\references


Abramowitz, M. and Stegun, I.A. (1968). \textit{Handbook of Mathematical Functions}, Dover, New York.

Bapat, R.B. (1989). Infinite divisibility of multivariate gamma distributions and \mbox{$M$--matrices,} \textit{Sankhy\={a}, Series A} \textbf{51}, 73--78.

Blacher, R. (2003).  Multivariate quadratic forms of random vectors, Journal of Multivariate Analysis \textbf{87},  2--23.

Fang, B.Q. (2005). Noncentral quadratic forms of the skew elliptical variables, \textit{Journal of Multivariate Analysis} \textbf{95}, 410--430.

Griffiths, R.C. (1984). Characterization of infinitely divisible multivariate gamma distributions, \textit{Journal of Multivariate Analysis} \textbf{15}, 13--20.

Khatri, C.G., Krishnaiah, P.R. and Sen, P.K. (1977). A note on the joint distribution of correlated quadratic forms, \textit{Journal of Statistical Planning and Inference} \textbf{1}, 299--307.

Krishnamoorthy, A.S. and Parthasarathy, M. (1951). A multivariate gamma type distribution, \textit{Annals of Mathematical Statistics} \textbf{22}, 549--557 (correction: ibid. (1960), \textbf{31}, p. 229).

Mathai, A.M. and Provost, S.B. (1992). \textit{Quadratic forms in random variables: Theory and applications}, Marcel Dekker, New York. 

Rao, C.R. (1973). \textit{Linear Statistical Inference and its Applications}, 2nd edition, Wiley, New York.

Royen, T. (1991). Expansions for the multivariate chi--square distribution, \textit{Journal of Multivariate Analysis} \textbf{38}, 213--232.

Royen, T. (1992).  On representation and computation of multivariate gamma distributions,
in: \textit{Data Analysis and Statistical Inference - Festschrift in Honour of Friedhelm Eicker}, 201--216,  Verlag Josef Eul, Bergisch Gladbach, K\"oln.

Royen, T. (1995). On some central and non--central multivariate chi--square distributions, \textit{Statistica Sinica} \textbf{5}, 373--397.

Royen, T. (1997). Multivariate gamma distributions (Update), \textit{Encyclopedia of Statistical Sciences}, Update Volume 1, 419--425, Wiley, New York.

Royen, T. (2006). Integral representations and approximations for multivariate gamma distributions, \textit{Annals of the Institute of Statistical Mathematics}, DOI 10.1007/s10463-006-0057-5.




\end{document}